\documentclass[11pt]{amsart}

\newtheorem{theorem}{Theorem}[section]

\newtheorem{lemma}[theorem]{Lemma}
\newtheorem{corollary}[theorem]{Corollary}

\title{A $5$-Engel associative algebra whose group of units is not $5$-Engel}

\author{Galina Deryabina}

\address{Department of Computational Mathematics and Mathematical Physics (FS-11), Bauman Moscow State Technical University, 2-nd Baumanskaya Street, 5, 105005 Moscow, Russia}

\author{Alexei Krasilnikov}

\address{Departamento de Matem\'atica, Universidade de
Bras\'\i lia, 70910-900 Bras\'\i lia, DF, Brazil}

\email{alexei@unb.br}

\date{}

\begin{document}

\maketitle

\begin{abstract}
Let $R$ be an associative ring with unity and let $[R]$ and $U(R)$ denote the associated Lie ring (with $[a,b]=ab-ba$) and the group of units of $R$, respectively. In 1983 Gupta and Levin proved that if $[R]$ is a nilpotent Lie ring of class $c$ then $U(R)$ is a nilpotent group of class at most $c$. The aim of the present note is to show that, in general, a similar statement does not hold if $[R]$ is $n$-Engel. We construct an algebra $R$ over a field of characteristic $\ne 2,3$ such that the Lie algebra $[R]$ is $5$-Engel  but the group $U(R)$ is not.
\end{abstract}

\section{Introduction}

Let $R$ be an associative ring with unity and let $[R]$ and $U(R)$ denote the associated Lie ring (with $[a,b]=ab-ba$) and the group of units of $R$, respectively. It is known that if $[R]$ is a nilpotent Lie ring of class $c$ then $U(R)$ is a nilpotent group of class at most $c$ (Gupta and Levin \cite{GL}). Also, if $[R]$ is metabelian then $U(R)$ is metabelian as well (Krasilnikov \cite{K} and Sharma and Srivastava \cite{SS}).

If $R$ is an associative ring such that $[R]$ is centre-by-metabelian then $U(R)$, in general, is not centre-by-metabelian. For instance, if $F$ is an infinite field of characteristic 2 and $R=M_2(F)$ is the algebra of all $2 \times 2$ matrices over $F$ then it is well-known that $[R]$ is centre-by-metabelian but $U(R)$ does not satisfy any non-trivial identity; in particular, $U(R)$ is not centre-by-metabelian.

However, suppose that $R$ is a unital associative algebra over a field of characteristic $0$ generated (as a unital algebra) by its nilpotent elements.  Then if $[R]$ is centre-by-metabelian then $U(R)$ is also centre-by metabelian (Krasilnikov and Riley \cite{KR}). Moreover, for such an algebra $R$, if $[R]$ satisfies an arbitrary multilinear Lie commutator identity then $U(R)$ satisfies the corresponding group commutator identity (see \cite{KR} for precise definitions); for example, if $[R]$ is solvable of length $n$ then $U(R)$ is also solvable of length at most $n$.

It is natural to ask whether a similar result holds for Lie commutator identities that are not multilinear; in particular, whether it holds for the Engel identity. We show that this is, in general, not the case.

More precisely, let $[x,y]=[x, _{(1)} y]=xy-yx$ and let $[x, _{(k+1)} y]=[[x, _{(k)}y],y]$
for $k \ge 1$. Recall that a Lie ring $L$ is $n$-Engel if
$[u, _{(n)} v]=0$ for all $u, v \in L$. Similarly, a group $G$ is
$n$-Engel if $(u,_{(n)} v)=1$ for all $u, v \in G$, where
$(x,y) = (x, _{(1)}y)= x^{-1} y^{-1}xy$ and
$(x, _{(k+1)} y) = ((x, _{(k)}y),y)$ for $k \ge 1$.
We are concerned with the following question.

\bigskip
\noindent
{\bf Question.}
\textit{Let $F$ be a field of characteristic $0$ and let $R$ be a unital associative $F$-algebra. Suppose that the Lie algebra $[R]$ is $n$-Engel. Is the group $U(R)$ also $n$-Engel?
}

\bigskip
If $n = 2$ then the answer to this question is ``yes''. Indeed, it is well-known (see, for instance, \cite[Theorem 3.1.1]{VL}) that if $[R]$ is $2$-Engel and $char F \ne 3$ then $[R]$ is nilpotent of class $2$. Hence, by \cite{GL}, $U(R)$ is nilpotent of class at most $2$ and, therefore, is $2$-Engel. If $n=3$ and the algebra $R$ is generated (as a unital algebra) by its nilpotent elements then the answer is ``yes'' as well; this can be deduced from the results of \cite{BRT} and \cite{Dickenschied} (see also \cite{AS}).

However, in general the answer to the question above is ``no''. Our result is as follows.

\begin{theorem}
\label{maintheorem1}
Let $F$ be a field of characteristic $\ne 2,3$. Then there is a unital associative $F$-algebra $R$ such that $[R]$ is a $5$-Engel Lie algebra but $U(R)$ is not a $5$-Engel group. This algebra $R$ is generated (as a unital $F$-algebra) by $2$ nilpotent elements.
\end{theorem}

Note that if $R$ is an associative unital ring and $[R]$ is $n$-Engel then $U(R)$ is $m$-Engel for some $m = m(n)$ (Riley and Wilson \cite{RW} and independently Amberg and Sysak \cite{AS}). If $R$ is an algebra over a field of characteristic $0$ then the existence of such $m = m(n)$ follows also from the results of Zelmanov \cite{Z} and Gupta and Levin \cite{GL}. One can check that in our example below the group $U(R)$ is $6$-Engel (and nilpotent of class $7$).

\medskip
We obtain Theorem \ref{maintheorem1} as a corollary of the result below about the adjoint groups of associative algebras.

Let $R$ be an associative ring with or without unity. It can be easily checked that $R$ is a monoid with respect to \textit{adjoint multiplication} defined by $u \circ v = u + v + uv$ $(u,v \in R)$. The group of units of this monoid is called the \textit{adjoint group} $R^{\circ}$ of $R$. It is well-known that if $R$ is nilpotent, that is, if $R^n = \{ 0 \}$ for some positive integer $n$, then $R^{\circ} = R$. On the other hand, if $R$ is a ring with unity $1$ then $R^{\circ}$ is isomorphic to the group of units $U(R)$ (the mapping $R^{\circ} \rightarrow R$ such that $a \rightarrow 1 + a$ is an isomorphism of $R^{\circ}$ onto $U(R)$).

Note that one can easily deduce from the results of \cite{GL,K,SS} that, for an associative ring $R$, if the Lie ring $[R]$ is nilpotent of class $c$ or metabelian then the adjoint group $R^\circ$ is also nilpotent of class at most $c$ or metabelian, respectively. Furthermore, if $R^\circ = R$ then the converse also holds: if $R^\circ$ is nilpotent of class $c$ then $[R]$ is nilpotent of class $c$ (Du \cite{D}) and if $R^\circ$ is metabelian then $[R]$ is metabelian (Amberg and Sysak \cite{AS3})

Let $F$ be a field and let $A$ be the free associative $F$-algebra without 1 on free generators $x,y$. Let $m(x,y), n(x,y) \in A$ be monic monomials in $x, y$. If $n(x,y) = m_1(x,y) m(x,y) m_2(x,y)$ for some monic monomials  $m_1(x,y), m_2(x,y) \in A \cup \{ 1 \}$ we say that $m(x,y)$ \textit{divides} $n(x,y)$ and $n(x,y)$ is a \textit{multiple} of $m(x,y)$.

Let $I$ be the ideal in $A$ generated by the following elements:

\begin{itemize}
\item[i)]
all monomials of degree 8;
\item[ii)]
all monomials of degree greater than  2 in $x$;
\item[iii)]
all monic monomials of degree 7 except $yxy^3xy$ and $y^2xyxy^2$;
\item[iv)]
all monic monomials of degree less than 7 which do not divide the
monomials $yxy^3xy$ and $y^2xyxy^2$;
\item[v)]
the polynomial
$2 xy^3xy - 5 yxyxy^2 - 2 yxy^3x + 5 y^2xyxy$;
\item[vi)]
the polynomial
$2 yxy^3xy - 5 y^2 x y x y^2$.
\end{itemize}

Let $B = A/I$. It can be easily seen that $B^8=0$. Thus, the associative algebra $B$ is nilpotent and, therefore,  $B^\circ = B$.

\begin{theorem}
\label{maintheorem2}
Let $F$ be a field of characteristic $ \ne 2,3$ and let $B=A/I$. Then $[B]$ is a $5$-Engel Lie algebra but the adjoint group $B^\circ$ is not $5$-Engel.
\end{theorem}

To deduce Theorem \ref{maintheorem1} from Theorem \ref{maintheorem2}, we embed a non-unital $F$-algebra $B$ into its unital hull. Let $B_1=F \oplus B$ be a direct sum of $F$-vector spaces $F$ and $B$. Then $B_1$ has a natural associative algebra structure in which the elements of $F$ act on $B$ by scalar multiplication. Further, the element $1$ of $F$ is unity for $B_1$ and the set $1+B$ forms a group under multiplication isomorphic to the adjoint group $B^\circ$ (the mapping $B^{\circ} \rightarrow B_1$ such that $u \rightarrow 1 + u$ is an isomorphism of $B^{\circ}$ onto $1+B$).

It can be easily checked that, since $[B]$ is $5$-Engel, so is the Lie ring $[B_1]$. On the other hand, since $B^\circ$ is not $5$-Engel, so are the subgroup $(1 + B) \simeq B^\circ$ of $U(B_1)$ and the group $U(B_1)$. Thus, we have

\begin{corollary}
\label{coroll}
Let $F$ be a field of characteristic $ \ne 2,3$ and let $B_1$ be the unital hull of $B$. Then the Lie algebra $[B_1]$ is $5$-Engel but the group $U(B_1)$ is not $5$-Engel.
\end{corollary}
Theorem \ref{maintheorem1} follows immediately from Corollary \ref{coroll} (with $R = B_1$).

\medskip
Recall that if $L$ is a nilpotent Lie algebra over a field of characteristic $0$ then $L$ is a group with the multiplication $*$ defined by the Baker-Campbell-Hausdorff formula: $x * y = log (e^x e^y) = x + y + \frac{1}{2}[x,y] + \dots$ (for details see, for example, \cite[\S 9.2]{Khukhro}). We denote this group by $L^*$. If $L^*$ is an $n$-Engel group for some $n \ge 1$ then $L$ is an $n$-Engel Lie algebra; this is well-known and can be deduced, for instance, from \cite[Lemma 10.12 (d)]{Khukhro}. However, the converse statement is false.

Indeed, if $B$ is a nilpotent associative algebra over a field of characteristic $0$ then it is well-known that $[B]^* \simeq B^\circ$ (see, for instance, \cite{KR}). Thus, by Theorem \ref{maintheorem2} we have

\begin{corollary}
Let $F$ be a field of characteristic $0$ and let $B = A/I$ be the associative $F$-algebra defined above. Let $L = [B]$. Then $L$ is a $5$-Engel Lie algebra such that the group $L^*$ is not $5$-Engel.
\end{corollary}

\noindent \textbf{Remarks.} 1. Theorems \ref{maintheorem1},  \ref{maintheorem2}, Corollary \ref{coroll} and their proofs remain valid for algebras over a unital associative and commutative ring $F$ such that $ 6 \ne 0$ in $F$.

2. One can check that if $R$ is a nilpotent associative algebra over an infinite field and the group $R^\circ$ is $n$-Engel then the Lie algebra $[R]$ is also $n$-Engel.

3. For an associative ring $R$, the Lie ring $[R]$ is nilpotent of class $c$ if and only if the adjoint semigroup $(R, \circ )$ is nilpotent of class $c$ in the sense of Mal'cev \cite{M} or Neumann-Taylor \cite{NT}. The ``only if'' part of this statement has been proved independently by Krasilnikov \cite{K2} and Riley and Tasic \cite{RT} and the ``if'' part by Amberg and Sysak \cite{AS2}. Note that if a group $G$, viewed as a semigroup, is Mal'cev or Neumann-Taylor nilpotent of class $c$ then $G$ is a nilpotent group of class $c$ in the usual sense \cite{M,NT}.

In \cite{R} Riley posed the following problem:

\medskip \noindent
\textit{Given any positive integer $n$, does there exist a semigroup variety $\mathcal P_n$ with the property that, for every associative ring $R$, the Lie ring $[R]$ is $n$-Engel if and only if the adjoint semigroup $(R, \circ )$ lies in $\mathcal P_n$?}

\medskip \noindent
This problem is not yet solved. Theorem \ref{maintheorem2} shows that, if such a variety of semigroups $\mathcal P_n$ exists, the groups that belong to $\mathcal P_n$ are not necessarily $n$-Engel.

\section{Proof of Theorem 1.2}

Let $I_0$ be the two-sided ideal in $A$ generated by the polynomials i)--iv) above. Let $C = A/I_0$. It is clear that $C = \oplus_{k=1}^{k=7} C_{(k)}$ where $C_{(k)}$ is the linear span of the monomials of degree $k$ in $x + I_0$, $y + I_0$.

It follows easily from the item iv) that
\begin{equation}
\label{eq_monom1}
x^2, \, xy^2x, \, y^4, \, y^2xy^2, \, xyxy^3, \, y^3xyx \in I_0.
\end{equation}
On the other hand, it can be easily checked that all monic monomials in $x, y$ satisfying iv) are multiples of the monomials above. Note that every monomial in $x, y$ of degree 1 in $x$ and of degree $5$ or $6$ in $y$ belongs to $I_0$ because such a monomial is a multiple of either $y^4$ or $y^2 x y^2$ and the latter monomials belong to $I_0$.  It is straightforward to check that an $F$-basis of $C_{(6)}$ is formed by the images of
 \[
 xy^3xy, \, yxyxy^2, \, yxy^3x, \, y^2xyxy
 \]
 and an $F$-basis of $C_{(7)}$ is formed by the images of
 \[
 yxy^3xy, \, y^2 x y x y^2.
 \]

By the definition of $I$, $I/I_0$ is the ideal of $C$ generated by the images of the polynomials
\[
h_1 = 2 xy^3xy - 5 yxyxy^2 - 2 yxy^3x + 5 y^2xyxy
\]
and
\[
h_2 = 2 yxy^3xy - 5 y^2 x y x y^2.
\]
Note that $x h_1 \equiv h_1 x \equiv 0 \pmod{I_0}$, $y h_1 \equiv - h_1 y \equiv h_2 \pmod{I_0}$ and $x h_2 \equiv h_2 x \equiv y h_2 \equiv h_2 y \equiv 0 \pmod{I_0}$. Hence, $h_1 + I_0$ and $h_2 + I_0$ form an $F$-basis of the ideal $I/I_0$. In particular, $C_{(7)} \cap I/I_0$ is a one-dimensional vector subspace in $C_{(7)}$ generated by the image of $h_2$ and $C_{(7)} /(C_{(7)} \cap I/I_0)$ is a one-dimensional vector space generated by the image of $y^2 x y x y^2$.

Let $a=x+I$, $b=y+I$. Then, by (\ref{eq_monom1}), we have
\begin{equation}
\label{eq_monom2}
a^2 = ab^2a = b^4 = b^2ab^2 = abab^3 = b^3aba = 0.
\end{equation}
It is clear that $B = \oplus_{k=1}^{k=7} B_{(k)}$ where $B_{(k)}$ is the linear span of monomials of degree $k$ in $a, b$. Note that $B_{(7)} \simeq C_{(7)} /(C_{(7)} \cap I/I_0)$ is a one-dimensional vector subspace in $B$ generated by $b^2 a b a b^2$.

To prove that $[B]$ is 5-Engel it suffices to check that $[u, {}_{(5)}v]=0$ for all $u, v \in B$. Let
\[
u = \alpha_1 a + \beta_1 b + \gamma_1 ab + \delta_1 ba + \mu_1 b^2 + u',
\qquad
v = \alpha_2 a + \beta_2 b + \gamma_2 ab + \delta_2 ba + \mu_2 b^2 + v',
\]
where $u'$ and $v'$ are linear combinations of monomials in $a, b$ of
degree at least 3. Then it is straightforward to check that
\begin{gather*}
[u, _{(5)}v] =
\alpha_1 [a, _{(5)} v] + \beta_1 [b, _{(5)} v]
+ \gamma_1 [ab, _{(5)} v] + \delta[ba, _{(5)} v] + \mu_1 [b^2, _{(5)} v]
\\
=
\alpha_1 \beta_2^5 f_0 + \alpha_1 \alpha_2 \beta_2^4 f_1
+ \alpha_1 \beta_2^4 \gamma_2 f_2 + \alpha_1 \beta_2^4 \delta_2 f_3
+ \alpha_1 \beta_2^4 \mu_2 f_4 + \alpha_1 \alpha_2 \beta_2^3 \mu_2 f_5
\\
+ \beta_1 \alpha_2 \beta_2^4 f_6 + \beta_1 \alpha_2^2 \beta_2^3 f_7
+ \beta_1 \beta_2^4 \gamma_2 f_8 + \beta_1 \alpha_2 \beta_2^3 \gamma_2 f_9
+ \beta_1 \beta_2^4 \delta_2 f_{10}
\\
+ \beta_1 \alpha_2 \beta_2^3 \delta_2 f_{11}
+ \beta_1 \alpha_2 \beta_2^3 \mu_2 f_{12}
+ \beta_1 \alpha_2^2 \beta_2^2 \mu_2 f_{13}
+ \gamma_1 \beta_2^5 f_{14}
\\
+ \gamma_1 \alpha_2 \beta_2^4 f_{15}
+ \delta_1 \beta_2^5 f_{16} + \delta_1 \alpha_2 \beta_2^4 f_{17}
+ \mu_1 \alpha_2 \beta_2^4 f_{18} + \mu_1 \alpha_2^2 \beta_2^3 f_{19},
\end{gather*}
where $f_i$ are multihomogeneous polynomials in $a, b$ of (total) degree $6$ or $7$.

Recall that every monomial in $x, y$ of degree 1 in $x$ and of degree $5$ or $6$ in $y$ belongs to $I_0$. Hence, every monomial in $a, b$ of degree 1 in $a$ and of degree 5 or 6 in $b$ is equal to 0. It follows immediately that $f_0 = f_4 = f_6 = f_8 = f_{10} = f_{12} = f_{14} = f_{16} = f_{18} = 0$.

It is straightforward to check that
\[
f_1 = [a, b, a, b, b, b] + [a, b, b, a, b, b] + [a, b, b, b, a, b]
+ [a, b, b, b, b, a],
\]
\begin{multline*}
f_2 = [a, ab, b, b, b, b] + [a, b, ab, b, b, b] + [a, b, b, ab, b, b]
\\
+ [a, b, b, b, ab, b] + [a, b, b, b, b, ab],
\end{multline*}
\begin{multline*}
f_3 = [a, ba, b, b, b, b] + [a, b, ba, b, b, b] + [a, b, b, ba, b, b]
\\
+ [a, b, b, b, ba, b] + [a, b, b, b, b, ba],
\end{multline*}
\begin{multline*}
f_5 = [a, b^2, a, b, b, b] + [a, b^2, b, a, b, b] + [a, b^2, b, b, a, b]
+ [a, b^2, b, b, b, a]
\\
+ [a, b, b^2, a, b, b] + [a, b, b^2, b, a, b]
+ [a, b, b^2, b, b, a] + [a, b, a, b^2, b, b]
\\
+ [a, b, b, b^2, a, b] + [a, b, b, b^2, b, a] + [a, b, a, b, b^2, b] + [a, b, b, a, b^2, b]
\\
+ [a, b, b, b, b^2, a] + [a, b, a, b, b, b^2] + [a, b, b, a, b, b^2]
+ [a, b, b, b, a, b^2],
\end{multline*}
\[
f_7 = [b, a, a, b, b, b] + [b, a, b, a, b, b] + [b, a, b, b, a, b]
+ [b, a, b, b, b, a] = - f_1,
\]
\begin{multline*}
f_9 = [b, a, ab, b, b, b] + [b, a, b, ab, b, b] + [b, a, b, b, ab, b]
+ [b, a, b, b, b, ab]
\\
+ [b, ab, a, b, b, b] + [b, ab, b, a, b, b]
+ [b, ab, b, b, a, b] + [b, ab, b, b, b, a]
\end{multline*}
\begin{multline*}
f_{11} = [b, a, ba, b, b, b] + [b, a, b, ba, b, b] + [b, a, b, b, ba, b]
+ [b, a, b, b, b, ba]
\\
+ [b, ba, a, b, b, b] + [b, ba, b, a, b, b]
+ [b, ba, b, b, a, b] + [b, ba, b, b, b, a],
\end{multline*}
\begin{multline*}
f_{13} = [b, a, b^2, a, b, b] + [b, a, b^2, b, a, b]
+ [b, a, b^2, b, b, a] + [b, a, a, b^2, b, b]
\\
+ [b, a, b, b^2, a, b]
+ [b, a, b, b^2, b, a] + [b, a, a, b, b^2, b] + [b, a, b, a, b^2, b]
\\
+ [b, a, b, b, b^2, a] + [b, a, a, b, b, b^2] + [b, a, b, a, b, b^2]
+ [b, a, b, b, a, b^2],
\end{multline*}
\begin{multline*}
f_{15} = [ab, a, b, b, b, b] + [ab, b, a, b, b, b] + [ab, b, b, a, b, b]
\\
+ [ab, b, b, b, a, b] + [ab, b, b, b, b, a],
\end{multline*}
\begin{multline*}
f_{17} = [ba, a, b, b, b, b] + [ba, b, a, b, b, b] + [ba, b, b, a, b, b]
\\
+ [ba, b, b, b, a, b] + [ba, b, b, b, b, a],
\end{multline*}
\[
f_{19} = [b^2, a, a, b, b, b] + [b^2, a, b, a, b, b]
+ [b^2, a, b, b, a, b] + [b^2, a, b, b, b, a].
\]

To proceed further we need the following lemma which is well-known and can
be easily proved by induction.

\begin{lemma}
\label{lem}
$ [x, {}_{(k)} y] = \sum_{i=0}^k {{k}\choose{i}} (-1)^i y^i x y^{(k-i)}$.
\end{lemma}

Now we will check that $f_1=0$. We have
$$
[a, b, a] = [ab-ba, a] = aba - b a^2 - a^2 b + aba = 2 aba
$$
because, by (\ref{eq_monom2}), $a^2 =0$.
Therefore,
$$
[a, b, a, b, b, b] = 2 [aba, b, b, b].
$$
By Lemma \ref{lem},
$$
[aba, b, b, b] = abab^3 - 3 babab^2 + 3 b^2abab - b^3aba,
$$
where, by (\ref{eq_monom2}), $abab^3=b^3aba=0$.
Hence,
$$
[aba, b, b, b] = - 3 babab^2 + 3 b^2abab
$$
and
\begin{equation}
\label{eq_comm11}
[a, b, a, b, b, b] = -6 babab^2 + 6 b^2abab.
\end{equation}

It is straightforward to check that $[a, b, b, a] = [a, b, a, b]$ so
\begin{equation}
\label{eq_comm12}
[a, b, b, a, b, b] = [a, b, a, b, b, b].
\end{equation}

Further, by Lemma \ref{lem},
$$
[a, b, b, b] = a b^3 - 3 b a b^2 + 3 b^2 a b - b^3 a
$$
so
$$
[a, b, b, b, a] = a b^3 a - 3 b a b^2 a + 3 b^2 a b a  - b^3 a^2
- a^2 b^3 + 3 a b a b^2 - 3 a b^2 a b + a b^3 a.
$$
Since, by (\ref{eq_monom2}),
$b a b^2 a = b^3 a^2 = a^2 b^3 = a b^2 a b =0$,
we have
$$
[a, b, b, b, a] = 2 a b^3 a + 3 b^2 a b a + 3 a b a b^2.
$$
Therefore,
$$
[a, b, b, b, a, b] = 2 a b^3 a b + 3 b^2 a b a b + 3 a b a b^3
- 2 b a b^3 a - 3 b^3 a b a - 3 b a b a b^2.
$$
By (\ref{eq_monom2}), $a b a b^3 = b^3 a b a = 0$ so
\begin{equation}
\label{eq_comm13}
[a, b, b, b, a, b]
= 2 a b^3 a b + 3 b^2 a b a b
- 2 b a b^3 a - 3 b a b a b^2.
\end{equation}
Finally, by Lemma \ref{lem} and (\ref{eq_monom2}),
$$
[a, b, b, b, b] = ab^4 - 4 b a b^3 + 6 b^2 a b^2 - 4 b^3 a b + b ^4 a
= -4 b a b^3 - 4 b^3 a b
$$
so, again by (\ref{eq_monom2}), we have
\begin{equation}
\label{eq_comm14}
[a, b, b, b, b, a]
= -4 b a b^3 a - 4  b^3 a b a + 4 a b a b^3  + 4 a b^3 a b
= -4 b a b^3 a + 4 a b^3 a b.
\end{equation}

Thus, by (\ref{eq_comm11}), (\ref{eq_comm12}), (\ref{eq_comm13}) and
(\ref{eq_comm14}), we have
$$
f_1 = -12 babab^2 + 12 b^2abab  + 2 a b^3 a b + 3 b^2 a b a b
- 2 b a b^3 a - 3 b a b a b^2 -4 b a b^3 a + 4 a b^3 a b
$$
$$
= -15 babab^2 + 15 b^2 a b a b  - 6 b a b^3 a + 6 a b^3 a b.
$$
By the item v) of the definition of the ideal $I$, $f_1 =0$.
Since $f_7 = - f_1$, we have $f_7 = 0$ as well.

One can check in a similar way using Lemma \ref{lem}, (\ref{eq_monom2})
and the item vi) of the definition of the ideal $I$ that
$f_2=f_3=f_5=f_9=f_{11}=f_{13}=f_{15}=f_{17}=f_{19}=0$.

More precisely, one can check using the relations (\ref{eq_monom2}) that
\[
[b^2, a, a, b, b, b] = [b^2, a, b, a, b, b] = [b^2, a, b, b, a, b] = [b^2, a, b, b, b, a] = 0.
\]
It follows that $f_{19}=0$. Similarly, $f_5 = f_{13} = 0$ because $f_5$ and $f_{13}$ are sums of certain commutators and one can check using (\ref{eq_monom2}) that all these commutators are equal to $0$.

Further, it can be checked using (\ref{eq_monom2}) that $f_2 = f_3 =0$ although the commutator summands of $f_2$ and $f_3$ are not, in general, equal to $0$. Finally, one needs the relations (\ref{eq_monom2}) as well as the item vi) of the definition of the ideal $I$ to check that $f_9 = f_{11} = f_{15} = f_{17} = 0$.

Thus, $[B]$ is a $5$-Engel Lie ring, as required.

\medskip
Now we prove that the group $B^\circ$ is not $5$-Engel. We will check that in $U(B_1)$
\[
( (1+a), _{(5)} (1+b)) = 1 + 6 b^2 a b a b^2 \ne 1.
\]
Hence, the subgroup $1+B$ of $U(B_1)$ is not $5$-Engel. Since $B^\circ$ is isomorphic to $1+B$, the group $B^\circ$ is not $5$-Engel as well.

Note that if $u \in B$ then $u^8 = 0$ so $(1 + u)^{-1} = 1 - u + u^2 - \dots - u^7$. It is straightforward to check that, for all $u, v \in B$,
\begin{equation}
\label{commutator1}
((1+u), (1+v)) = 1 + [u,v] - u^2 v + uvu + v^2 u - vuv + w
\end{equation}
where $w$ is a linear combination of monomials of degree at least $4$ in $u, v$. It follows that if $u = u(a,b)$ is a linear combination of some monomials of degree $k \ge 2$ in $a, b$ and, possibly, some monomials of degree $>k$ then
\begin{equation}
\label{commutator2}
((1+u), (1+b)) = 1 + [u,b] + b^2 u - bub + w' = 1 + [u,b] - b [u,b] + w'
\end{equation}
where $w' \in B^{k+3}$.

By (\ref{commutator1}), we have
\begin{multline*}
((1+a), (1 + b)) = 1 + [a,b] + aba - a^2 b + b^2 a - bab + w_1
\\
=  1 + [a,b] + aba + b^2 a - bab + w_1
\end{multline*}
where $w_1 \in B^4$.

Let $u_1 = [a,b] + aba + b^2 a - bab + w_1$; then $((1+a), (1 + b)) = 1 + u_1$. By (\ref{commutator2}), we have
\begin{multline*}
((1+a), _{(2)}(1 + b)) = ((1 + u_1), (1+b)) = 1 + [u_1, b] - b [u_1, b] + w_2'
\\
= 1 + [a,b,b] + [(aba + b^2 a - bab), b] - b [a,b, b] +  w_2
\\
= 1 + [a,b,b] +abab - baba - 2 b [a,b,b] + w_2
\end{multline*}
where $w_2', w_2 \in B^5$.

Similarly, one can check that
\begin{multline*}
((1+a), _{(3)}(1 + b)) = 1 + [a,b,b,b] + abab^2 - 2 babab
\\
+ b^2aba - 3 b [a,b,b,b] + w_3
\end{multline*}
where $w_3 \in B^6$,
\begin{multline*}
((1+a), _{(4)}(1 + b)) = 1 + [a,b,b,b,b] - 3 babab^2 + 3 b^2abab - 4 b [a,b,b,b,b] + w_4
\end{multline*}
where $w_4 \in B^7$ and
\begin{multline*}
((1+a), _{(5)}(1 + b)) = 1 + [a,b,b,b,b,b] + 6 b^2abab^2 - 5 b [a,b,b,b,b,b] + w_5
\end{multline*}
where $w_5 = 0$ because $w_5 \in B^8$ and $B^8 = 0$ and $[a,b,b,b,b,b] = 0$ because the Lie algebra $[B]$ is $5$-Engel.

Thus, $( (1+a), _{(5)} (1+b)) = 1 + 6 b^2 a b a b^2 $. Since $B_{(7)}$ is a one-dimensional vector subspace in $B$ generated by $b^2 a b a b^2$, we have $b^2 a b a b^2 \ne 0$ and $( (1+a), _{(5)} (1+b)) \ne 1$. It follows that $B^\circ$ is not $5$-Engel group, as required.

The proof of Theorem \ref{maintheorem2} is completed.

\end{document}